\documentclass[12pt]{article}
\usepackage[centertags]{amsmath}
\usepackage{amsfonts}
\usepackage{amssymb}
\usepackage{graphicx}
\newtheorem{theorem}{Theorem}[section]

\newtheorem{definition}[theorem]{Definition}

\usepackage[verbose,colorlinks=true,naturalnames=true,linkcolor=blue,]{hyperref}

\newcounter{rown}

\begin{document}

\title{Twisted Quantum Deformations of \\Lorentz and Poincar\'{e} algebras\footnote{Invited
talk at the VII International Workshop ''Lie Theory and its Applications in Physics'',
18--24 June 2007, Varna, Bulgaria}
}
\author{V.N. Tolstoy
\footnote{Supported by the grants RFBR-05-01-01086 and FNRA NT05-241455GIPM
.}\\
\\ Institute of Nuclear Physics, Moscow State University, \\
119 992 Moscow, Russia; e-mail: tolstoy@nucl-th.sinp.msu.ru}

\date{}

\maketitle
\begin{abstract}
We discussed twisted quantum deformations of $D=4$ Lorentz and Poincar\'{e} algebras. In
the case of Poincar\'{e} algebra it is shown that almost all classical $r$-matrices of
S.~Zakrzewski classification can be presented as a sum of subordinated $r$-matrices of
Abelian and Jordanian types. Corresponding twists describing  quantum deformations are
obtained in explicit form. This work is an extended version of the paper
\url{arXiv:0704.0081v1 [math.QA]}.

\end{abstract}

\section{Introduction}
The quantum deformations of relativistic symmetries are described by Hopf-algebraic
deformations of Lorentz and  Poincar\'{e} algebras. Such quantum deformations are
classified by Lorentz and Poincar\'{e} Poisson structures. These Poisson structures given
by classical $r$-matrices were classified already some time ago by S. Zakrzewski in
\cite{Z1} for the Lorentz algebra and in \cite{Z2} for the Poincar\'{e} algebra. In the
case of the Lorentz algebra a complete list of classical $r$-matrices involves the four
independent formulas and the corresponding quantum deformations in different forms were
already discussed in literature (see \cite{ChD, M, BLT1, BLT2, BLT3}). In the case of
Poincar\'{e} algebra the total list of the classical $r$-matrices, which satisfy the
homogeneous classical Yang-Baxter equation, consists of 20 cases which have various
numbers of free parameters. Analysis of these twenty solutions shows that eighteen 
of them can be presented as a sum of subordinated $r$-matrices which  
are of Abelian and Jordanian types. Corresponding twists describing  quantum deformations
are obtained in explicit form.

This work is extended version of the paper \cite{T1}.

\setcounter{equation}{0}
\section{Preliminaries}

Let $r$ be a classical $r$-matrix of a Lie algebra $\mathfrak{g}$, i.e.
$r\in\,\stackrel{2}\wedge\mathfrak{g}$ and $r$ satisfies to the classical Yang--Baxter
equation (CYBE)
\begin{eqnarray}\label{p1}
[r^{12},\,r^{13}+\,r^{23}] + [r^{13},\,r^{23}]\!\!&=\!\!&\Omega~,
\end{eqnarray}
where $\Omega$ is $\mathfrak{g}$-invariant element, $\Omega\in(\stackrel{3}\wedge
\mathfrak{g})_{\mathfrak{g}}$. We consider two types of the classical $r$-matrices and
corresponding twists.

Let the classical $r$-matrix $r=r_{A}^{}$ has the form
\begin{eqnarray}\label{p2}
r_{A}^{}\!\!&=\!\!&\sum_{i=1}^{n}y_{i}\wedge x_{i}~,
\end{eqnarray}
where all elements $x_i, y_i$ $(i=1,\ldots,n)$ commute among themselves. Such an
$r$-matrix is called of Abelian type. The corresponding twist is given as follows
\begin{eqnarray}\label{p3}
F_{r_{A}^{}}\!\!&=\!\!&\exp\frac{\tilde{r}_{A}^{}}{2}=
\exp\Bigl(\frac{1}{2}\sum_{i=1}^{n}x_{i}\wedge y_{i}\Bigr)~.
\end{eqnarray}
This twisting two-tensor $F:=F_{r_{A}^{}}$ satisfies the cocycle equation
\begin{equation}\label{p4}
F^{12}(\Delta\otimes{\rm id})(F)\;=\;F^{23}({\rm id}\otimes\Delta)(F)~,
\end{equation}
and the "unital" normalization condition
\begin{equation}\label{p5}
(\epsilon \otimes{\rm id})(F)\;=\;({\rm id}\otimes\epsilon )(F)\;=\;1~.
\end{equation}
The twisting element $F$ defines a deformation of the universal enveloping algebra
$U(\mathfrak{g})$ considered as a Hopf algebra. The new deformed coproduct and antipode
are given as follows
\begin{equation}\label{p6}
\Delta^{(F)}(a)\;=\;F\Delta(a)F^{-1}~,\qquad S^{(F)}(a)=uS(a)u^{-1}
\end{equation}
for any $a\in U(\mathfrak{g})$, where $\Delta(a)$ is a co-product before twisting, and
\begin{eqnarray}\label{p7}
u\!\!&=\!\!&\sum_i f^{(1)}_{i}S(f^{(2)}_i)
\end{eqnarray}
if $F=\sum_i f^{(1)}_i\otimes f^{(2)}_i$.

Let the classical $r$-matrix $r=r_{J_n}^{}(\xi)$ has the form\footnote{Here entering the
parameter deformation $\xi$ is a matter of convenience.}
\begin{equation}\label{p8}
r_{J_n}^{}(\xi)\;=\;\xi\,\Bigl(\sum_{\nu=0}^{n}y_{\nu}\wedge x_{\nu}\Bigr)~,
\end{equation}
where the elements $x_\nu,y_\nu$ $(\nu=0,1,\ldots, n)$ satisfy the relations\footnote{It
is easy to verify that the two-tensor (\ref{p8}) indeed satisfies the homogenous
classical Yang-Baxter equation (\ref{p1}) (with $\Omega=0$), if the elements
$x_\nu,y_\nu$ $(\nu=0,1,\ldots,n)$ are subject to the relations (\ref{p9}).}
\begin{equation}\label{p9}
\begin{array}{rcl}
[x_{0},y_{0}]\!\!\!&=\!\!\!&y_{0},\quad\;\;
[x_{0},x_{i}]\,=\,(1-t_{i})x_{i},\quad\;\;[x_{0},y_{i}]\,=\,t_{i}y_{i},
\\[7pt]
[x_{i},y_{j}]\!\!\!&=\!\!\!&\delta_{ij}y_{0},\;\;[x_{i},x_{j}]\,=\,
[y_{i},y_{j}]\,=\,0,\;\;[y_{0},x_{j}]\,=\,[y_{0},y_{j}]=0,
\end{array}
\end{equation}
$(i,j=1,\ldots,n)$, $(t_{i}\in{\mathbb C})$. Such an $r$-matrix is called of Jordanian
type. The corresponding twist is given as follows \cite{T2,T3}
\begin{eqnarray}\label{p10}
F_{r_{J_n}^{}}\!\!&=\!\!& \exp\Big(\xi\sum\limits_{i=1}^{n}x_{i}\otimes
y_{i}\;e^{-2t_{i}\sigma}\Bigr)\exp(2x_{0}^{}\otimes\sigma)~,
\end{eqnarray}
where $\sigma:=\frac{1}{2}\ln(1+\xi y_{0})$.\footnote{Similar twists for Lie algebras
$\mathfrak{sl}(n)$, $\mathfrak{so}(n)$ and $\mathfrak{sp}(2n)$ were firstly constructed
in the papers \cite{KLM, KLO, LSK, AKL}.}

{\it Remark.} The zero component $r_{J_0}^{}(\xi):=\xi y_{0}\wedge x_{0}$ in (\ref{p8})
itself is the classical Jordanian $r$-matrix and the corresponding Jordanian twist is
given by the formula (\ref{p10}) for $n=0$, i.e.
$F_{r_{J_0}^{}}=\,\exp(2x_{0}^{}\otimes\sigma)$.

Let $r$ be an arbitrary $r$-matrix of $\mathfrak{g}$. We denote a support of $r$ by
$\mathop{\rm Sup}(r)$\footnote{The support $\mathop{\rm Sup}(r)$ is a subalgebra of
$\mathfrak{g}$ generated by the elements $\{x_i,y_i\}$ if $r=\sum_{i}y_i\wedge x_i$.}.
The following definition is useful.
\begin{definition}
Let $r^{}_1$ and $r^{}_2$ be two arbitrary classical $r$-matrices. We say that $r^{}_2$
is subordinated to $r^{}_1$, $r^{}_1\succ r^{}_2$, if $\delta_{r^{}_1}(\mathop{\rm
Sup}(r^{}_2))=0$, i.e.
\begin{equation}\label{p11}
\delta_{r_{1}^{}}(x)\;:=\;[x\otimes1+1\otimes x,\,r_{1}^{}]\;=\;0~, \quad \forall x\in
\mathop{\rm Sup}(r_2)~.
\end{equation}
\end{definition}
If $r^{}_1\succ r^{}_2$ then $r=r^{}_1+r^{}_2$ is also a classical $r$-matrix (see
\cite{BD}). The subordination enables us to construct a correct sequence of
quantizations. For instance, if the $r$-matrix of Jordanian type (\ref{p8}) is
subordinated to the $r$-matrix of Abelian type (\ref{p2}), $r_{\!A}^{}\succ r_{\!J}^{}$,
then the total twist corresponding to the resulting $r$-matrix $r=r_{\!A}^{}+r_{\!J}^{}$
is given as follows
\begin{eqnarray}\label{p12}
F_{r}\!\!&=\!\!&F_{r_{\!J}^{}}\,F_{r_{\!A}^{}}.
\end{eqnarray}
The further definition is also useful.
\begin{definition}
A twisting two-tensor $F_{r}(\xi)$ of a Hopf algebra $U(\mathfrak{g})$, satisfying the
conditions (\ref{p4}) and (\ref{p5}), is called locally $r$-symmetric if the expansion of
$F_{r}(\xi)$ in powers of the parameter deformation $\xi$ has the form
\begin{eqnarray}\label{p13}
F_{r}(\xi)\!\!&=\!\!& 1+c\,r+\mathcal{O}(\xi^2)\ldots
\end{eqnarray}
where $r$ is a classical $r$-matrix, and $c$ is a numerical coefficient, $c\neq0$.
\end{definition}
It is evident that the Abelian twist (\ref{p3}) is globally $r$-symmetric and the twist
of Jordanian type (\ref{p10}) does not satisfy the relation (\ref{p13}), i.e. it is not
locally $r$-symmetric.  It is a matter of direct verification to prove the following
theorem.
\begin{theorem}\label{t13}
Let $F_{r}(\xi)$ be a twisting two-tensor of a Hopf algebra $U(\mathfrak{g})$, and
$\omega:=\sqrt{u}$ where $u$ is give by the formula (\ref{p7}), then the new twisting
two-tensor $F_{r}^{(\omega)}(\xi)$,
\begin{eqnarray}\label{p14}
F_{r}^{(\omega)}(\xi)\!\!&:=\!\!&\omega^{-1}\otimes\omega^{-1} F_{r}(\xi)\Delta(\omega)~,
\end{eqnarray}
is locally $r$-symmetric.
\end{theorem}

\setcounter{equation}{0}
\section{Quantum deformations of Lorentz algebra}

The results of this section in different forms were already discussed in literature (see
\cite{ChD, M, BLT1, BLT2}).

The classical canonical basis of the $D=4$ Lorentz algebra, $\mathfrak{o}(3,1)$, can be
described by the anti-Hermitian six generators ($h$, $e_{\pm}$, $h'$, $e'_{\pm}$)
satisfying the following non-vanishing commutation relations\footnote{Since the real Lie
algebra $\mathfrak{o}(3,1)$ is standard realification of the complex Lie
$\mathfrak{sl}(2,\mathbb{C})$ these relations are easy obtained from the defining
relations for $\mathfrak{sl}(2,\mathbb{C})$, i.e. from (\ref{L1}).}:
\begin{eqnarray}
&[h,\,e_{\pm}^{}]\;=\;\pm e_{\pm}^{}\,,\qquad [e_{+}^{},\,e_{-}^{}]\;=\;2h~,
\label{L1}
\\[5pt]
&[h,\,e'_{\pm}]\;=\;\pm e'_{\pm}~,\quad\;\; [h',\,e_{\pm}]\;=\;\pm e'_{\pm}~,\quad\;\;
[e_{\pm}^{},\,e'_{\mp}]\;=\;\pm2h'~, \label{L2}
\\[5pt]
& [h',\,e'_{\pm}]\;=\;\mp e_{\pm}^{}~,\qquad [e'_{+},\,e'_{-}]\;=\;-2h~,
\label{L3}
\end{eqnarray}
and moreover
\begin{equation}\label{L4}
x^*\;=\;-x\qquad (\forall\;x\;\in\; \mathfrak{o}(3,1))~.\textbf{}
\end{equation}
A complete list of classical $r$-matrices which describe all Poisson structures and
generate quantum deformations for $\mathfrak{o}(3,1)$ involve the four independent
formulas \cite{Z1}:
\begin{eqnarray}\label{L5}
r_{1}^{}\!\!&=\!\!&\alpha e_{+}\wedge h~,
\\[5pt]\label{L6}
r_{2}^{}\!\!&=\!\!&\alpha(e_{+}\wedge h-e'_{+}\wedge h')+2\beta e'_{+}\wedge e_{+}~,
\\[5pt]\label{L7}
r_{3}^{}\!\!&=\!\!&\alpha(e'_{+}\wedge e_{-}\!+e_{+}\wedge e'_{-})+\beta(e_{+}\wedge
e_{-}\!-e'_{+}\wedge e'_{-})-2\gamma h\wedge h',
\\[5pt]\label{L8}
r_{4}^{}\!\!&=\!\!&\alpha\bigl(e'_{+}\wedge e_{-}\!+e_{+}\wedge e'_{-}\!-2h\wedge
h'\bigr)\pm e_{+}\wedge e'_{+}~.
\end{eqnarray}
If the universal $R$-matrices of the quantum deformations corresponding to the classical
$r$-matrices (\ref{L5})--(\ref{L8}) are unitary then these $r$-matrices are
anti-Hermitian, i.e.
\begin{equation}\label{L9}
r^*_{j}\;=\;-r_{j}\qquad (j=1,2,3,4)~.
\end{equation}
Therefore the $*$-operation (\ref{L4}) should be lifted to the tensor product
$\mathfrak{o}(3,1)\otimes \mathfrak{o}(3,1)$. There are two variants of this lifting:
{\it direct} and {\it flipped}, namely,
\begin{eqnarray}\label{L10}
(x\otimes y)^*\!\!&=\!\!&x^*\otimes y^*\qquad({\rm *-direct})~,
\\[5pt]\label{L11}
(x\otimes y)^*\!\!&=\!\!&y^*\otimes x^*\qquad({\rm *-flipped})~.
\end{eqnarray}
We see that if the "direct" lifting of the $*$-operation (\ref{L4}) is used then all
parameters in (\ref{L5})--(\ref{L8}) are pure imaginary. In the case of the "flipped"
lifting (\ref{L11}) all parameters in (\ref{L5})--(\ref{L8}) are real.

The first two $r$-matrices (\ref{L5}) and (\ref{L6}) satisfy the homogeneous CYBE and
they are of Jordanian type. If we assume (\ref{L10}), the corresponding quantum
deformations were described detailed in the paper \cite{BLT2} and they are entire defined
by the twist of Jordanian type:
\begin{equation}\label{L12}
F_{r_1^{}}^{}\,=\,\exp{(2h\otimes\sigma)}~,\quad\sigma\,=\,\frac{1}{2}\ln(1+\alpha e_{+})
\end{equation}
for the $r$-matrix (\ref{L5}), and
\begin{equation}\label{L13}
F_{r_2^{}}^{}\,=\,\exp{\Bigl(\frac{\imath\beta}{\alpha^2}\;\sigma\wedge\varphi\Bigr)}\,
\exp{(h\otimes\sigma-h'\otimes\varphi)}~,
\end{equation}
\begin{equation}\label{L14}
\sigma \,=\,\frac{1}{2}\ln\left[(1+\alpha e_+)^2\!+(\alpha e'_+)^2
\right],\quad\varphi\,=\,\arctan{\frac{\alpha e'_+}{1+\alpha e_+}}
\end{equation}
for the $r$-matrix (\ref{L6}). It should be recalled that the twists (\ref{L12}) and
(\ref{L13}) are not locally $r$-symmetric. Using the representation of the Jordanian
twist $F_{r_1^{}}^{}$ in binomial series form (see \cite{KST1, KST2})
\begin{equation}\label{L15}
\begin{array}{rcl}
F_{r_1^{}}^{}\!\!&=\!\!&(1\otimes1+1\otimes\alpha e_+)^{\frac{h}{2}\otimes 1}
\\[10pt]
\!\!&=\!\!&1+\! \sum_{k>0}\frac{\alpha^k}{k!}h(h-1)\cdots(h-k+1)\otimes e_{+}^{k}
\end{array}
\end{equation}
we can easy obtain the explicit form of the element (\ref{p7}):
\begin{equation}\label{L16}
u\,=\,1+\!\sum_{k>0}\frac{(-\alpha)^k}{k!}h(h-1)\cdots(h-k+1)e_{+}^{k}=\exp{(-\alpha
he_+)}
\end{equation}
and therefore
\begin{equation}\label{L17}
\omega:=\sqrt{u}\,=\,\exp{\Bigl(-\frac{1}{2}\alpha he_+\Bigr)}~.
\end{equation}
By Theorem \ref{t13} the locally $r$-symmetric Jordanian twist is given as
follows\footnote{Another form of the locally $r$-symmetric Jordanian twist was presented
in \cite{Oh} (see also \cite{T1}).}
\begin{equation}\label{L18}
\begin{array}{rcl}
F_{r_1^{}}^{(\omega)}\!\!&=\!\!&\displaystyle\exp{\Bigl(\frac{\alpha}{2}\bigr(he_+\otimes1+1\otimes
he_+\bigr)\Bigr)}\, \exp{(2h\otimes\sigma)}\times
\\[10pt]
&&\times\displaystyle\exp{\Bigl(-\frac{\alpha}{2}\bigl(he_+\otimes1+h\otimes
e_{+}+e_{+}\otimes h+1\otimes he_+\bigr)\Bigr)}~.
\end{array}
\end{equation}
In a similar way one can obtain the locally $r$-symmetric expression for the twist
$F_{r_2^{}}^{}$, (\ref{L13}).

The last two $r$-matrices (\ref{L7}) and (\ref{L8}) satisfy the non-homogeneous
(modified) CYBE and they can be easily obtained from the solutions of the complex algebra
$\mathfrak{o}(4,\mathbb{C})\simeq \mathfrak{sl}
(2,\mathbb{C})\oplus\mathfrak{sl}(2,\mathbb{C})$ which describes the  complexification of
$\mathfrak{o}(3,1)$. Indeed, let us introduce the complex basis of Lorentz algebra
$(\mathfrak{o}(3,1)\simeq\mathfrak{sl}(2;\mathbb{C})\oplus\mathfrak
{\overline{sl}}(2,\mathbb{C}))$ described by two commuting sets of complex generators:
\begin{eqnarray}\label{L19}
H_1\!\!&=\!\!& \frac{1}{2}\,(h+\imath h')~,\qquad
E_{1\pm}\;=\;\frac{1}{{2}}\,(e_{\pm}^{}+ \imath e'_{\pm})~,
\\[5pt]\label{L20}
H_2\!\!&=\!\!& \frac{1}{{2}}\, (h-\imath h')~,\qquad
E_{2\pm}\;=\;\frac{1}{{2}}\,(e_{\pm}^{}-\imath e'_{\pm})~,
\end{eqnarray}
which satisfy the non-vanishing relations (compare with (\ref{L1}))
\begin{equation}\label{L21}
[H_k,\,E_{k\pm}]\;=\;\pm E_{k\pm}~,\qquad [E_{k+},\,E_{k-}]\;=\;2 H_k\qquad(k=1,2)~.
\end{equation}
The $*$-operation describing the real structure acts on the generators $H_k$, and
$E_{k\pm}$ ($k=1,2$) as follows
\begin{equation}\label{L22}
H_1^*=-H_2^{}~,\quad E_{1\pm}^*=- E_{2\pm}^{}~,\quad H_2^*=-H_1^{}~,\quad E_{2\pm}^*=
-E_{1\pm}^{}~.
\end{equation}
The classical $r$-matrix $r_3$, (\ref{L7}), and $r_4$, (\ref{L8}), in terms of the
complex basis (\ref{L19}), (\ref{L20}) take the form
\begin{equation}\label{L23}
\begin{array}{rcl}
r_3^{}\!\!&=\!\!&r'_3+r''_3~,
\\[7pt]
r'_3\!\!&:=\!\!&2(\beta+\imath\alpha)E_{1+}\wedge
E_{1-}+2(\beta-\imath\alpha)E_{2+}\wedge E_{2-}~,
\\[7pt]
r''_3\!\!&:=\!\!&4\imath\gamma\,H_{2}\wedge H_{1}~,
\end{array}
\end{equation}
and
\begin{equation}\label{L24}
\begin{array}{rcl}
r_4^{}\!\!&=\!\!&r'_4+r''_4~,\quad
\\[7pt]
r'_4\!\!&:=\!\!&2\imath\alpha(E_{1+}\wedge E_{1-}-E_{2+}\wedge E_{2-}-2H_{1}\wedge
H_{2})~,
\\[7pt]
r''_4\!\!&:=\!\!&2\imath\lambda\,E_{1+}\wedge E_{2+}. 
\end{array}
\end{equation}

For the sake of convenience we introduce parameter\footnote{We can reduce this parameter
$\lambda$ to $\pm1$ by automorphism of $\mathfrak{o}(4,\mathbb{C})$.}$\lambda$ in
$r''_{4}$. It should be noted that $r'_{3}$, $r''_{3}$ and $r'_{4}$, $r''_{4}$ are
themselves classical $r$-matrices. We see that the $r$-matrix $r'_{3}$ is simply a sum of
two standard $r$-matrices of $\mathfrak{sl}(2;\mathbb{C})$, satisfying the anti-Hermitian
condition $r^*=-r$. Analogously, it is not hard to see that the $r$-matrix $r_4$
corresponds to a Belavin-Drinfeld triple \cite{BD} for the Lie algebra
$\mathfrak{sl}(2;\mathbb{C})\oplus \mathfrak{\overline{sl}}(2,\mathbb{C}))$. Indeed,
applying the Cartan automorphism $E_{2\pm}\rightarrow E_{2\mp}$, $H_{2}\rightarrow-H_{2}$
we see that this is really correct (see also \cite{IO2001}).

We firstly describe quantum deformation corresponding to the classical $r$-matrix
$r_{3}$, (\ref{L23}). Since the $r$-matrix $r''_{3}$ is Abelian and it is subordinated to
$r'_{3}$ therefore the algebra $\mathfrak{o}(3,1)$ is firstly quantized in the direction
$r'_{3}$ and then an Abelian twist corresponding to the $r$-matrix $r''_{3}$ is applied.
We introduce the complex notations $z_{\pm}:=\beta\pm\imath\alpha$. It should be noted
that $z_{-}^{}=z_{+}^*$ if the parameters $\alpha$ and $\beta$ are real, and
$z_{-}^{}=-z_{+}^*$ if the parameters $\alpha$ and $\beta$ are pure imaginary. From
structure of the classical $r$-matrix $r_{3}'$ it follows that a quantum deformation
$U_{r'_{1}} (\mathfrak{o}(3,1))$ is a combination of two $q$-analogs of
$U(\mathfrak{sl}(2;\mathbb{C}))$ with the parameter $q_{z_{+}^{}}$ and $q_{z_{-}^{}}$,
where $q_{z_{\pm}^{}}:=\exp{z_{\pm}^{}}$. Thus $U_{r'_{3}}(\mathfrak{o}(3,1))\cong
U_{q_{z_{+}^{}}^{}}(\mathfrak{sl}(2;\mathbb{C}))\otimes U_{q_{z_{-}^{}}^{}}
(\overline{\mathfrak{sl}}(2;\mathbb{C}))$ and the standard generators $q_{z_{+}^{}}^{\pm
H_{1}}$, $E_{1\pm}$ and $q_{z_{-}^{}}^{\pm H_{2}}$, $E_{2\pm}$ satisfy the following
non-vanishing defining relations
\begin{eqnarray}\label{L25}
q_{z_{+}}^{H_1}E_{1\pm}\!\!\!&=\!\!\!&q_{z_{+}}^{\pm1}E_{1\pm}\,q_{z_{+}}^{H_1},\quad
\;\; [E_{1+},\,E_{1-}]\;=\;\frac{q_{z_{+}}^{2H_1}-q_{z_{+}}^{-2H_1}}
{q_{z_{+}}^{}-q_{z_{+}}^{-1}},
\\[7pt]\label{L26}
q_{z_{-}}^{H_2}E_{2\pm}\!\!&=\!\!&q_{z_{-}}^{\pm1}E_{2\pm}\,q_{z_{-}}^{H_2},\quad\;\;
[E_{2+},\,E_{2-}]\;=\;\frac{q_{z_{-}}^{2H_2}-q_{z_{-}}^{-2H_2}}
{q_{z_{-}}^{}-q_{z_{-}}^{-1}}~.
\end{eqnarray}
In this case the co-product $\Delta_{r'_{1}}$ and antipode $S_{r'_{1}}$ for the
generators $q_{z_{+}^{}}^{\pm H_{1}}$, $E_{1\pm}$ and $q_{z_{-}^{}}^{\pm H_{2}}$,
$E_{2\pm}$ can be given by the formulas:
\begin{eqnarray}\label{L27}
\Delta_{r'_{1}}^{}(q_{z_{+}}^{\pm H_{1}})\!\!&=\!\!&q_{z_{+}}^{\pm H_{1}}\otimes
q_{z_{+}}^{\pm H_{1}},
\\[7pt]\label{L28}
\Delta_{r'_{1}}^{}(E_{1\pm})\!\!&=\!\!&E_{1\pm}\otimes q_{z_{+}}^{H_{1}}+
q_{z_{+}}^{-H_{1}}\otimes E_{1\pm},
\\[7pt]\label{L29}
\Delta_{r'_{1}}^{}(q_{z_{-}}^{\pm H_{2}})\!\!\!&=\!\!\!&q_{z_{-}}^{\pm H_{2}}\otimes
q_{z_{-}}^{\pm H_{2}},
\\[7pt]\label{L30}
\Delta_{r'_{1}}^{}(E_{2\pm})\!\!&=\!\!&E_{2\pm}\otimes q_{z_{-}}^{H_{2}}+
q_{z_{-}}^{-H_{2}}\otimes E_{2\pm},\phantom{aaaaa}
\end{eqnarray}
\begin{eqnarray}\label{L31}
S_{r'_{1}}^{}(q_{z_{+}}^{\pm H_{1}})\!\!&=\!\!&q_{z_{+}}^{\mp H_{1}},\quad
S_{r'_{1}}^{}(E_{1\pm})\;=\;-q_{z_{+}}^{\pm1}E_{1\pm}~,
\\[7pt]\label{L32}
S_{r'_{1}}^{}(q_{z_{-}^{}}^{\pm H_{2}})\!\!&=\!\!&q_{z_{-}}^{\mp H_{2}},\quad
S_{r'_{1}}^{}(E_{2\pm})\;=\;-q_{z_{-}}^{\,\pm1}E_{2\pm}~.
\end{eqnarray}
The $*$-involution describing the real structure on the generators (\ref{L22}) can be
adapted to the quantum generators as follows
\begin{equation}\label{L33}
\begin{array}{rcccl}
(q_{z_{+}}^{\pm H_{1}})^*\!\!&=\!\!&q_{z_{+}^{*}}^{\mp
H_{2}},\qquad\;E_{1\pm}^*\!\!&=\!\!&- E_{2\pm}^{}~,
\\[7pt]
(q_{z_{-}}^{\pm H_{2}})^*\!\!&=\!\!&q_{z_{-}^{*}}^{\mp H_{1}},\qquad
E_{2\pm}^*\!\!&=\!\!&- E_{1\pm}^{}~,
\end{array}
\end{equation}
and there exit two $*$-liftings: {\it direct} and {\it flipped}, namely,
\begin{eqnarray}\label{L34}
(a\otimes b)^*\!\!&=\!\!&a^*\otimes b^*\qquad({\rm*-direct})~,
\\[5pt]\label{L35}
(a\otimes b)^*\!\!&=\!\!&b^*\otimes a^*\qquad({\rm*-flipped})
\end{eqnarray}
for any $a\otimes b\in U_{r'_{3}}(\mathfrak{o}(3,1))\otimes U_{r'_{3}}
(\mathfrak{o}(3,1))$, where the $*$-direct involution corresponds to the case of the pure
imaginary parameters $\alpha,\,\beta$ and the $*$-flipped involution corresponds to the
case of the real deformation parameters $\alpha,\,\beta$.  It should be stressed that the
Hopf structure on $U_{r'_{3}}(\mathfrak{o}(3,1))$ satisfy the consistency conditions
under the $*$-involution
\begin{equation}\label{L36}
\Delta_{r'_{3}}(a^*)=(\Delta_{r'_{3}}(a))^*,\;\; S_{r'_{3}}((S_{r'_{3}}
(a^*))^{*})=a\;\;(\forall x\in U_{r'_{3}}(\mathfrak{o}(3,1)).
\end{equation}

Now we consider deformation of the quantum algebra $U_{r'_{3}}(\mathfrak{o}(3,1))$
(secondary quantization of $U(\mathfrak{o}(3,1))$) corresponding to the additional
$r$-matrix $r''_{3}$, (\ref{L23}). Since the generators $H_{1}$ and $H_{2}$ have the
trivial coproduct
\begin{eqnarray}\label{L37}
\Delta_{r'_{3}}(H_{k})\!\!&=\!\!&H_{k}\otimes 1+1\otimes H_{k}\quad(k=1,2)~,
\end{eqnarray}
therefore the unitary two-tensor
\begin{eqnarray}\label{L38}
F_{r_{3}''}^{}\!\!:=\!\!&q_{\imath\gamma}^{H_{1}\wedge H_{2}}\qquad
(F_{r_{3}''}^*\;=\;F_{r_{1}''}^{-1})
\end{eqnarray}
satisfies the cocycle condition (\ref{p4}) and the "unital" normalization condition
(\ref{p5}). Thus the complete deformation corresponding to the $r$-matrix $r_{3}^{}$ is
the twisted deformation of $U_{r'_{3}}(\mathfrak{o}(3,1))$, i.e. the resulting coproduct
$\Delta_{r_{3}^{}}$ is given as follows
\begin{eqnarray}\label{L39}
\Delta_{r_{3}^{}}^{}(x)\!\!&=\!\!&F_{r_{1}''}^{}\Delta_{r_{1}'}^{}(x)
F_{r_{3}''}^{-1}\quad(\forall x\in U_{r'_{1}}(\mathfrak{o}(3,1))~.
\end{eqnarray}
and in this case the resulting antipode $S_{r_{3}^{}}^{}$ does not change,
$S_{r_{3}^{}}^{}=S_{r'_{3}}^{}$. Applying the twisting two-tensor (\ref{L38}) to the
formulas (\ref{L27})--(\ref{L30}) we obtain
\begin{eqnarray}\label{L40}
\Delta_{r_{3}^{}}(q_{z_{+}}^{\pm H_{1}})\!\!&=\!\!&q_{z_{+}}^{\pm H_{1}}\otimes
q_{z_{+}}^{\pm H_{1}},
\\[7pt]\label{L41}
\Delta_{r_{3}^{}}(E_{1+})\!\!&=\!\!&E_{1+}\otimes q_{z_{+}}^{H_{1}}
q_{\imath\gamma}^{H_{2}}+q_{z_{+}}^{-H_{1}} q_{\imath\gamma}^{-H_{2}}\otimes E_{1+}~,
\\[7pt]\label{L42}
\Delta_{r_{3}^{}}(E_{1-})\!\!&=\!\!&E_{1-}\otimes q_{z_{+}}^{H_{1}}
q_{\imath\gamma}^{-H_{2}}+q_{z_{+}}^{-H_{1}}q_{\imath\gamma}^{H_{2}}\otimes E_{1-}~,
\\[7pt]\label{L43}
\Delta_{r'_{1}}(q_{z_{-}}^{\pm H_{2}})\!\!&=\!\!&q_{z_{-}}^{\pm H_{2}}\otimes
q_{z_{-}}^{\pm H_{2}},
\\[7pt]\label{L44}
\Delta_{r_{3}}(E_{2+})\!\!&=\!\!&E_{2+}\otimes q_{z_{-}}^{H_{2}}
q_{\imath\gamma}^{-H_{1}}+q_{z_{-}}^{-H_{2}} q_{\imath\gamma}^{H_{1}}\otimes E_{2+}~,
\\[7pt]\label{L45}
\Delta_{r_{3}}(E_{2-})\!\!&=\!\!&E_{2-}\otimes q_{z_{-}}^{H_{2}}
q_{\imath\gamma}^{H_{1}}+q_{z_{-}}^{-H_{2}}q_{\imath\gamma}^{-H_{1}}\otimes E_{2-}~.
\end{eqnarray}

Next, we describe quantum deformation corresponding to the classical $r$-matrix $r_{4}$
(\ref{L24}). Since the $r$-matrix $r_{4}'(\alpha):=r_4'$ is a particular case of
$r_{3}^{}(\alpha,\beta,\gamma):=r_{3}^{}$, namely
$r_{4}'(\alpha)=r_{3}^{}(\alpha,\beta=0,\gamma=\alpha)$, therefore a quantum deformation
corresponding to the $r$-matrix $r_4'$ is obtained from the previous case by setting
$\beta=0,\gamma=\alpha$, and we have the following formulas for the coproducts
$\Delta_{r_{4}'}$:
\begin{eqnarray}\label{L46}
\Delta_{r_{4}'}(q_{\xi}^{\pm H_{k}})\!\!&=\!\!&q_{\xi}^{\pm H_{k}}\otimes q_{\xi}^{\pm
H_{k}} \qquad (k=1,2)~,
\\[7pt]\label{L47}
\Delta_{r_{4}'}(E_{1+})\!\!&=\!\!&E_{1+}\otimes q_{\xi}^{H_{1}+H_{2}}+
q_{\xi}^{-H_{1}-H_{2}}\otimes E_{1+}~,
\\[7pt]\label{L48}
\Delta_{r_{4}'}(E_{1-})\!\!&=\!\!&E_{1-}\otimes q_{\xi}^{H_{1}-H_{2}}+
q_{\xi}^{-H_{1}+H_{2}}\otimes E_{1-}~,
\\[7pt]\label{L49}
\Delta_{r_{4}'}(E_{2+})\!\!&=\!\!&E_{2+}\otimes q_{\xi}^{-H_{1}-H_{2}}+
q_{\xi}^{H_{1}+H_{2}}\otimes E_{2+}~,
\\[7pt]\label{L50}
\Delta_{r_{4}'}(E_{2-})\!\!&=\!\!&E_{2-}\otimes q_{\xi}^{H_{1}-H_{2}}+
q_{\xi}^{-H_{1}+H_{2}}\otimes E_{2-}~,
\end{eqnarray}
where we set $\xi:=\imath\alpha$.

Consider the two-tensor
\begin{eqnarray}\label{L51}
F_{r_{4}''}\!\!:=\!\!&\exp_{q_{\xi}^{2}}^{}\big(-2\imath\lambda
E_{1+}q_{\xi}^{H_{1}+H_{2}}\otimes E_{2+}q_{\xi}^{H_{1}+H_{2}}\big)~.
\end{eqnarray}
Using properties of $q$-exponentials (see \cite{KT1}) is not hard to verify that
$F_{r_{4}''}$ satisfies the cocycle equation (\ref{p4}). Thus the quantization
corresponding to the $r$-matrix $r_4$ is the twisted $q$-deformation
$U_{r_{4}'}(\mathfrak{o}(3,1))$. Explicit formulas of the co-products
$\Delta_{r_{4}^{}}^{}(\cdot)=F_{r_{4}''}^{}\Delta_{r_{4}'}^{}(\cdot) F_{r_{4}''}^{-1}$
and antipodes $S_{r_4}(\cdot)$ will be presented in the outgoing paper \cite{BLT3}.

 \setcounter{equation}{0}
\section{Quantum deformations of Poincar\'{e} algebra}

The Poincar\'{e} algebra ${\mathcal{P}}(3,1)$ of the 4-dimensional space-time is
generated by 10 elements: the six-dimensinal Lorentz algebra
$\mathfrak{o}(3,1)$ with the generators  $M_i$, $N_i$ ($i=1,2,3$):
\begin{equation}\label{P1}
\begin{array}{rcl}
[M_i,\,M_j ]\!\!&=\!\!&\imath\epsilon_{ijk}\,M_k~,
\\[7pt]
[M_i,\,N_j]\!\!&=\!\!&\imath\epsilon_{ijk}\,N_k~,
\\[7pt]
[N_i,\,N_j]\!\!&=\!\!&-\imath\epsilon_{ijk}\,M_k~,
\end{array}
\end{equation}
and the four-momenta $P_0$, $P_j$ $(j=1,2,3)$ with the standard commutation relations:
\begin{eqnarray}\label{P2}
[M_j,\,P_k]\!\!&=\!\!&\imath\epsilon_{jkl}\,P_l~,\qquad [M_j,\,P_0]\;=\;0~,
\\[5pt]\label{P3}
[N_j,\,P_k]\!\!&=\!\!&-\imath\delta_{jk}\,P_0~,\quad\;\;[N_j,\,P_0]\;=\;-\imath P_j^{}~.
\end{eqnarray}
The physical generators of the Lorentz algebra, $M_i$, $N_i$ ($i = 1,2,3$), are related
with the canonical basis $h,h',e_{\pm},e'_{\pm}$ as follows
\begin{eqnarray}\label{P4}
h\!\!&=\!\!&\imath N_3~,\qquad e_{\pm}\;=\;\imath (N_1\pm\,M_2),
\\[5pt]\label{P5}
h'\!\!&=\!\!&\imath M_3~,\qquad e'_{\pm}\;=\;\imath (M_1\mp N_2).
\end{eqnarray}
The subalgebra generated by the four-momenta $P_0$, $P_j$ $(j=1,2,3)$ will be denoted by
$\mathbf{P}$, and we also set $P_{\pm}:=P_{0}\pm P_{3}$.

S.~Zakrzewski has shown in \cite{Z2} that each classical $r$-matrix,
$r\in\mathcal{P}(3,1) \wedge\mathcal{P}(3,1)$,  has a decomposition
\begin{equation}\label{P6}
r=a+b+c~,
\end{equation}
where $a\in\mathbf{P}\wedge\mathbf{P}$, $b\in\mathbf{P}\wedge\mathfrak{o}(3,1)$,
$c\in\mathfrak{o}(3,1)\wedge\mathfrak{o}(3,1)$ satisfy the following relations
\begin{eqnarray}\label{P7}
[[c,c]]\!\!&=\!\!&0~,
\\[3pt]\label{P8}
[[b,c]]\!\!&=\!\!&0~,
\\[3pt]\label{P9}
2[[a,c]]+[[b,b]]\!\!&=\!\!&t\Omega\quad (t\in \mathbb{R},\;\Omega\neq0)~,
\\[3pt]\label{P10}
[[a,b]]\!\!&=\!\!&0~.
\end{eqnarray}
Here $[[\cdot,\cdot]]$ means the Schouten bracket. Moreover a total list of the classical
$r$-matrices for the case $c\neq0$ and also for the case $c=0$, $t=0$ was
found.\footnote{Classification of the $r$-matrices for the case $c=0$, $t\neq0$ is an
open problem up to now.} The results are presented in the following table taken from
\cite{Z2}:
\\[10pt]
{
\begin{tabular}{ccccc}
\hline $c$ & $b$ & $a$ & $\#$ & $N$\\
\hline $\gamma h'\wedge h$ & $0$ & $\alpha P_{+}\wedge P_{-}+\tilde{\alpha}P_{1}\wedge
P_{2}$ & $2$ & $1$\\
\hline $\gamma e'_{+}\wedge e_{+}$ & $\beta_{1}b_{P_{+}}^{}+\beta_{2}P_{+}\wedge h'$ &
$0$ & $1$ & $2$\\
$$ & $\beta_{1} b_{P_{+}}^{}$ & $\alpha P_{+}\wedge P_{1}$ & $1$ & $3$\\
$$ & $\gamma\beta_{1}(P_{1}\wedge e_{+}+P_{2}\wedge e'_{+})$ & $P_{+}\wedge(\alpha_{1}P_{1}\!+
\alpha_{2}P_{2})-\gamma\beta_{1}^2P_{1}\wedge P_{2}$ & $2$ & $4$\\
\hline $\gamma(h\wedge e_{+}$ & $$ & $$ & $$ & $$\\
$-h'\wedge e'_{+})$ & $0$ & $0$ & $1$ & $5$\\
$+\gamma_{1}e'_{+}\wedge e_{+}$ & $$ & $$ & $$ & $$\\
\hline $\gamma h\wedge e_{+}$ & $\beta_{1}b_{P_{2}}^{}+\beta_{2}P_{2}\wedge e_{+}$ & $0$
& $1$ & $6$\\
\hline $0$ & $\beta_{1}b_{P_{+}}^{}+\beta_{2}P_{+}\wedge h'$ & $0$ & $1$ & $7$\\
$$ & $\beta_{1}b_{P_{+}}^{}+\beta_{2}P_{+}\wedge e_{+}$ & $0$ & $1$ & $8$\\
$$ & $P_{1}\wedge(\beta_{1}e_{+}+\beta_{2}e'_{+})\,+$ & $\alpha P_{+}\wedge P_{2}$ &
$2$ & $9$\\
$$ & $\beta_{1}P_{+}\wedge(h+\chi e_{+}),\;\chi=0,\pm1$ & $$ & $$ & $$\\
$$ & $\beta_{1}(P_{1}\wedge e'_{+}+P_{+}\wedge e_{+})$ & $\alpha_{1}P_{-}\wedge
P_{1}+\alpha_{2} P_{+}\wedge P_{2}$ & $2$ & $10$\\
$$ & $\beta_{1} P_{2}\wedge e_{+}$ & $\alpha_{1}P_{+}\wedge P_{1}+\alpha_{2}P_{-}\wedge P_{2}$ &
$1$ & $11$\\
$$ & $\beta_{1} P_{+}\wedge e_{+}$ & $P_{-}\!\wedge(\alpha P_{+}\!+\!\alpha_{1}P_{1}\!+\!\alpha_{2}P_{2})\!+
\tilde{\alpha} P_{+}\!\wedge P_{2}$ & $3$ & $12$\\
$$ & $\beta_{1} P_{0}\wedge h'$ & $\alpha_{1}P_{0}\wedge P_{3}+\alpha_{2}P_{1}\wedge P_{2}$ &
$2$ & $13$\\
$$ & $\beta_{1} P_{3}\wedge h'$ & $\alpha_{1}P_{0}\wedge P_{3}+\alpha_{2}P_{1}\wedge P_{2}$ &
$2$ & $14$\\
$$ & $\beta_{1} P_{+}\wedge h'$ & $\alpha_{1}P_{0}\wedge P_{3}+\alpha_{2}P_{1}\wedge P_{2}$ &
$1$ & $15$\\
$$ & $\beta_{1} P_{1}\wedge h$ & $\alpha_{1}P_{0}\wedge P_{3}+\alpha_{2}P_{1}\wedge P_{2}$ &
$2$ & $16$\\
$$ & $\beta_{1} P_{+}\wedge h$ & $\alpha P_{1}\wedge P_{2}+\alpha_{1}P_{+}\wedge P_{1}$ &
$1$ & $17$\\
$$ & $P_{+}\wedge(\beta_{1} h+\beta_{2} h')$ & $\alpha_{1} P_{1}\wedge P_{2}$ & $1$ & $18$\\
\cline{2-5}
$$ & $0$ & $\alpha_{1} P_{1}\wedge P_{+}$ & $0$ & $19$\\
$$ & $$ & $\alpha_{1} P_{1}\wedge P_{2}$ & $0$ & $20$\\
$$ & $$ & $\alpha_{1} P_{0}\wedge P_{3}+\alpha_{2}P_{1}\wedge P_{2}$ & $1$ & $21$\\
\hline
\end{tabular}
}
\begin{center}
{\bf Table 1.} Normal forms of $r$ for $c\neq0$ or $t=0$.
\end{center}
Where $b_{P_{+}}^{}$ and $b_{P_{2}}^{}$ are given as follows:
\begin{eqnarray}\label{P11}
b_{P_{+}}^{}\!\!&=\!\!&P_{1}\wedge e_{+}-P_{2}\wedge e'_{+}+P_{+}\wedge h~,
\\[5pt]\label{P12}
b_{P_{2}}^{}\!\!&=\!\!&2P_{1}\wedge h'+P_{-}\wedge e'_{+}-P_{+}\wedge e'_{-}~.
\end{eqnarray}
The table lists 21 cases labelled by the number $N$ in the last column. In the forth
column (labelled by $\#$) the number of essential parameters (more precisely - the
maximal number of such parameter) involved in deformation are indicated. This number is
in any cases less than the number of parameters actually using in the table. Moreover we
introduced an additional parameter $\gamma$ in the component $a$ (in the cases 2,3,4,5,6)
and a parameter $\beta_1^{}$ in the component $b$ (in the cases 7--18) and also a
parameter $\alpha_1^{}$ in the component $a$ (in the cases 19--21)\footnote{In the
original paper by S. Zakrzewski \cite{Z2} all these additional parameters are equal to 1
and the numbers in the forth column correspond to this situation.}. The final reduction
of the number of the actual parameters can be achieved using of automorphisms of the
Poincar\'{e} algebra ${\mathcal{P}}(3,1)$ (see details in \cite{Z2}).

Now we would like to show that each $r$-matrices of the given table, $r_{N}$ ($1\leq
N\leq 21$), can be presented as a sum of subordinated $r$-matrices which almost all are
of Abelian and Jordanian types and we write down possible twisting two-tensors. Firstly
we analyze classical $r$-matrices for the case $c\neq0$.

1). The first $r$-matrix $r_1$,
\begin{eqnarray}\label{P13}
r_{1}^{}\!\!&=\!\!&\gamma h'\wedge h+\alpha P_{+}\wedge P_{-}+\tilde{\alpha} P_{1}\wedge
P_{2}~,
\end{eqnarray}
is a sum of two subordinated Abelian $r$-matrices
\begin{eqnarray}\label{P14}
\begin{array}{rcl}
r_1\!\!&=\!\!&r_1'+r_1''~,\quad r_1'\,\succ\,r_1''~,
\\[7pt]
r_1'\!\!&=\!\!&\alpha\,P_{+}\wedge P_{-}+\tilde{\alpha}\,P_{1}\wedge P_{2}~,
\\[7pt]
r_1''\!\!&=\!\!&\gamma h'\wedge h~.
\end{array}
\end{eqnarray}
Therefore the total twist defining quantization in the direction to this $r$-matrix is
the ordered product of  two the Abelian twits
\begin{eqnarray}\label{P15}
F_{r_1^{}}\!\!&=\!\!&F_{r_1''}\,F_{r_1'}\;=\;\exp\bigl(\gamma h\wedge h'\bigr)
\exp\bigl(\alpha P_{-}\wedge P_{+}+\tilde{\alpha}P_{2}\wedge P_{1}\bigr)~.
\end{eqnarray}

2). The second  $r$-matrix $r_2$,
\begin{equation}\label{P16}
\begin{array}{rcl}
r_{2}^{}\!\!&=\!\!&\gamma e'_{+}\wedge e_{+}+\beta_{1}(P_{1}\wedge e_{+}-P_{2}\wedge
e'_{+}+P_{+}\wedge h)\,+
\\[7pt]
&& +\,\beta_{2}P_{+}\wedge h'~,
\end{array}
\end{equation}
is a sum of three subordinated $r$-matrices where one of them is of Jordanian type and
two are of Abelian type
\begin{eqnarray}\label{P17}
\begin{array}{rcl}
r_{2}^{}\!\!&=\!\!&r_{2}'+r_{2}''+r_{2}'''~,\quad r_{2}'\,\succ\,r_{2}''~,\quad
r_{2}'+r_{2}''\succ\, r_{2}'''~,
\\[7pt]
r_{2}'\!\!&=\!\!&\beta_{1}(P_{1}\wedge e_{+}-P_{2}\wedge e'_{+}+P_{+}\wedge h)~,
\\[7pt]
r_{2}''\!\!&=\!\!&\gamma e'_{+}\wedge e_{+}~,
\\[7pt]
r_{2}'''\!\!&=\!\!&\beta_{2}P_{+}\wedge h'~.
\end{array}
\end{eqnarray}
The corresponding twisting two-tensor is given by the following formula
\begin{eqnarray}\label{P18}
F_{r_{2}^{}}\!\!&=\!\!&F_{r_{2}'''}\,F_{r_{2}''}\,F_{r_{2}'}~,
\end{eqnarray}
where
\begin{eqnarray}\label{P19}
\begin{array}{rcl}
F_{r_{2}^{'}}\!\!&=\!\!& \exp\bigr(\beta_{1}(e_{+}^{}\otimes P_{1}-e_{+}'\otimes
P_{2})\bigr)\exp(2h\otimes\sigma_{+})~,
\\[12pt]
F_{r_{2}^{''}}\!\!&=\!\!&\exp(\gamma e_{+}\wedge e'_{+})~,
\\[7pt]
F_{r_{2}'''}\!\!&=\!\!&\displaystyle\exp(\frac{\beta_2}{\beta_1}h'\wedge \sigma_{+})~.
\end{array}
\end{eqnarray}
Here and below we set $\sigma_+:=\frac{1}{2}\ln(1+\beta_{1}P_{+})$.

3). The third $r$-matrix $r_3$,
\begin{eqnarray}\label{P20}
\begin{array}{rcl}
r_3\!\!&=\!\!&\gamma e'_{+}\wedge e_{+}+\beta_{1}(P_{1}\wedge e_{+}-P_{2}\wedge e'_{+}
+P_{+}\wedge h)\,+
\\[7pt]
&&+\,\alpha P_{+}\wedge P_{1}~,
\end{array}
\end{eqnarray}
can be presented as a sum of three subordinated $r$-matrices where one of them is of
Jordanian type and two are of Abelian type
\begin{eqnarray}\label{P21}
\begin{array}{rcl}
r_{3}^{}\!\!&=\!\!&r_{3}'+r_{3}''+r_{3}'''~,\quad r_{3}'\;\succ\,
r_{3}''~,\quad~,r_{3}'+r_{3}''\;\succ\,r_{3}'''~,
\\[9pt]
r_{3}'\!\!&=\!\!&\displaystyle P_{1}\wedge(\beta_{1}e_{+}-\alpha P_{+})-\beta_{1}P_{2}
\wedge e'_{+}+\beta_{1}P_{+}\wedge h~,
\\[5pt]
r_{3}''\!\!&=\!\!&\displaystyle\gamma e'_{+}\wedge\bigl(e_{+}-
\frac{\alpha}{\beta_{1}}P_{+}\bigr)~,
\\[7pt]
r_{3}'''\!\!&=\!\!&\displaystyle\frac{\gamma\alpha}{\beta_{1}}e'_{+}\wedge P_{+}~.
\end{array}
\end{eqnarray}
The corresponding twist is given by the following formula
\begin{eqnarray}\label{P22}
F_{r_3^{}}\!\!&=\!\!&F_{r_{3}''}\,F_{r_{3}'}~,
\end{eqnarray}
where
\begin{eqnarray}\label{P23}
\begin{array}{rcl}
F_{r_{3}^{'}}\!\!&=\!\!& \exp\bigr((\beta_{1} e_{+}^{}-\alpha P_{+})\otimes
P_{1}-\beta_{1} e'_{+}\otimes P_{2}\bigr)\exp(2h\otimes\sigma_{+})~,
\\[9pt]
F_{r_{3}^{''}}\!\!&=\!\!&\displaystyle\exp\Bigl(\gamma\bigl(e_{+}-
\frac{\alpha}{\beta_{1}}P_{+}\bigr)\wedge e'_{+}\Bigr)~,
\\[9pt]
F_{r_{3}^{''}}\!\!&=\!\!&\displaystyle\exp\Bigl(
\frac{\gamma\alpha}{\beta_{1}^2}\sigma_{+}\wedge e'_{+}\Bigr)~.
\end{array}
\end{eqnarray}

4). The fourth $r$-matrix $r_4$,
\begin{eqnarray}\label{P24}
\begin{array}{rcl}
r_{4}^{}\!\!&=\!\!&\gamma(e'_{+}\wedge e_{+}+\beta_{1} P_{1}\wedge e_{+}+\beta_{1}
P_{2}\wedge e'_{+}- \beta_{1}^{2}P_{1}\wedge P_{2})\,+
\\[7pt]
&&+\,P_{+}\wedge(\alpha_{1}P_{1}+\alpha_{2}P_{2})~,
\end{array}
\end{eqnarray}
is a sum of two subordinated $r$-matrices of Abelian type
\begin{eqnarray}\label{P25}
\begin{array}{rcl}
r_{4}^{}\!\!&=\!\!&r_{4}'+r_{4}''~,\quad r_{4}'\;\succ\, r_{4}''~,
\\[7pt]
r_{4}'\!\!&=\!\!&P_{+}\wedge(\alpha_{1}P_{1}+\alpha_{2}P_{2})~,
\\[7pt]
r_{4}''\!\!&=\!\!&\gamma(e'_{+}+\beta_{1} P_{1})\wedge(e_{+}^{}-\beta_{1} P_{2})~.
\end{array}
\end{eqnarray}
The corresponding twist is given by the following formula
\begin{eqnarray}\label{P26}
F_{r_4^{}}\!\!&=\!\!&F_{r_{4}''}\,F_{r_{4}'}~,
\end{eqnarray}
where
\begin{eqnarray}\label{P27}
\begin{array}{rcl}
F_{r_{4}^{'}}\!\!&=\!\!& \exp\bigr((\alpha_{1}P_{1}+\alpha_{2}P_{2})\wedge
P_{+}^{}\bigr)~,
\\[9pt]
F_{r_{4}^{''}}\!\!&=\!\!& \exp\bigr(\gamma(e_{+}-\beta_{1}P_{1})\wedge(e'_{+}+\beta_{1}
P_{2})\bigr)~.
\end{array}
\end{eqnarray}

{\it Remark.} The parameter $\beta_{1}$ can be removed by the rescaling automorphism
$\beta_{1}P_{\nu}\rightarrow P_{\nu}$ ($\nu=0,1,2,3$).

5). The fifth $r$-matrix $r_5$ is the $r$-matrix of the Lorentz algebra, (\ref{L6}) of
the Section~3, and the corresponding twist is given by the formulas (\ref{L13}) and
(\ref{L14}) of the previous Section 3.

6). The sixth  $r$-matrix $r_6$,
\begin{eqnarray}\label{P28}
\begin{array}{rcl}
r_{6}^{}\!\!&=\!\!&\gamma h\wedge e_{+}+\beta_{1}(2P_{1}\wedge h'+P_{-}\wedge
e'_{+}-P_{+}\wedge e'_{-})\,+
\\[7pt]
&&+\,\beta_{2}P_{2}\wedge e_{+}~,
\end{array}
\end{eqnarray}
is a sum of three subordinated $r$-matrices
\begin{eqnarray}\label{P29}
\begin{array}{rcl}
r_6\!\!&=\!\!&r_{6}'+r_{6}''+r_{6}''',\quad r_{6}'\,\succ\,r_{6}''~,\quad
r_{6}'+r_{6}''\,\succ\,r_{6}'''~,
\\[7pt]
r_{6}'\!\!&=\!\!&\beta_{1}(2P_{1}\wedge h'+P_{-}\wedge e'_{+}-P_{+}\wedge e'_{-})~,
\\[7pt]
r_{6}''\!\!&=\!\!&\gamma h\wedge e_{+}~,
\\[7pt]
r_{6}'''\!\!&=\!\!&\beta_{2}P_{2}\wedge e_{+}~.
\end{array}
\end{eqnarray}
The $r$-matrix $r_{6}''$ is of Jordanian type and the $r$-matrix $r_{6}'''$ is of Abelian
type while the first $r$-matrix $r_{6}'$ is not Abelian and Jordanian type because the
total $r$-matrix (\ref{P28}) satisfy the non-homogeneous classical Yang-Baxter equation
(\ref{P9}) with $t\neq0$. In terms of the generators  $M_i$, $N_i$ ($i=1,2,3$) the
$r$-matrix $r_{6}'$ has the form
\begin{eqnarray}\label{P30}
r'_{6}\!\!&=\!\!&2\imath\beta_{1}(P_{1}\wedge M_{3}-P_{3}\wedge M_{1}-P_{0}\wedge
N_{2})~.
\end{eqnarray}
Unfortunately a quantum deformation corresponding to this $r$-matrix is unknown to us,
and it is very likely that it can be obtain by contraction procedure from the $q$-analog
$U_q(\mathfrak{so}(5))$ in the same way as for the $\kappa$-Poincar'e deformation (see
\cite{LNRT}).

Now we will analyze classical $r$-matrices when $c=0,\;t=0$. In this case there are
twelve solutions with $b\neq0$, $r_{N}$ ($7\leq N\leq 18$),  and three solutions with
$b=0$, $r_{N}$ ($19\leq N\leq 21$).

7). The seventh  $r$-matrix $r_7$,
\begin{eqnarray}\label{P31}
r_{7}^{}\!\!&=\!\!&\beta_{1}(P_{1}\wedge e_{+}-P_{2}\wedge e'_{+}+ P_{+}\wedge
h)+\beta_{2}P_{+}\wedge h'~,
\end{eqnarray}
is a partial case of the $r$-matrix $r_{2}(\gamma,\beta_1,\beta_2):=r_{2}$, namely
$r_{7}=r_{2}(\gamma=0,\beta_1,\beta_2)$, therefore the corresponding twist is defined by
the formulas (\ref{P18}) and (\ref{P19}) where $\gamma=0$.

8). The eighth $r$-matrix $r_8$,
\begin{eqnarray}\label{P32}
r_{8}^{}\!\!&=\!\!&\beta_{1}(P_{1}\wedge e_{+}-P_{2}\wedge e'_{+}+ P_{+}\wedge
h)+\beta_{2}P_{+}\wedge e_{+}~,
\end{eqnarray}
is a sum of two subordinated $r$-matrices of Jordanian and Abelian types
\begin{eqnarray}\label{P33}
\begin{array}{rcl}
r_{8}^{}\!\!&=\!\!&r_{8}'+r_{8}''~,\quad r_{8}'\,\succ\, r_{8}''~,
\\[7pt]
r_{8}'\!\!&=\!\!&\beta_{1}(P_{1}\wedge e_{+}-P_{2}\wedge e'_{+}+P_{+}\wedge h)~,
\\[7pt]
r_{8}''\!\!&=\!\!&\beta_{2}P_{+}\wedge e_{+}~.
\end{array}
\end{eqnarray}
The corresponding twisting two-tensor is given by the following formula
\begin{eqnarray}\label{P34}
F_{r_{8}^{}}\!\!&=\!\!&F_{r_{8}''}\,F_{r_{8}'}~,
\end{eqnarray}
where
\begin{eqnarray}\label{P35}
\begin{array}{rcl}
F_{r_{8}'}\!\!&=\!\!&\exp\bigr(\beta_{1}(e_{+}^{}\otimes P_{1}-e_{+}'\otimes
P_{2})\bigr)\exp(2h\otimes\sigma_{+})~,
\\[9pt]
F_{r_{8}''}\!\!&=\!\!&\displaystyle\exp(\frac{\beta_{2}}{\beta_{1}}e_{+}
\wedge\sigma_{+})~.
\end{array}
\end{eqnarray}

9). The two-tensor $\tilde{r}_9^{}$, which corresponds to the ninth row of {\bf Table 1},
has the form:
\begin{eqnarray}\label{P36}
\begin{array}{rcl}
\tilde{r}_{9}^{}\!\!&=\!\!&P_{1}\wedge(\beta_{1}e_{+}+\beta_{2}e'_{+})+\beta_{1}P_{+}
\wedge(h+\chi e_{+})
\\[7pt]
\!\!&\!\!&+\,\alpha P_{+}\wedge P_{2}~,
\end{array}
\end{eqnarray}
Unfortunately we could not prove that $\tilde{r}_{9}^{}$ for $\chi\neq0$ satisfies the
equation (\ref{P9}) and thus we consider the modified variant $r_{9}^{}$:
\begin{eqnarray}\label{P37}
\begin{array}{rcl}
r_{9}^{}\!\!&=\!\!&P_{1}\wedge(\beta_{1}e_{+}+\beta_{2}e'_{+})+P_{+}\wedge\bigl
(\beta_{1}h+ \chi(\beta_{1}e_{+}+\beta_{2}e'_{+})\bigr)
\\[7pt]
\!\!&\!\!&+\,\alpha P_{+}\wedge P_{2}~,
\end{array}
\end{eqnarray}
This two-tensor satisfies the equations (\ref{P9}) and (\ref{P10}) and therefore it is a
classical $r$-matrix. We can present $r_{9}^{}$ as a sum of two subordinated $r$-matrices
of Jordanian and Abelian types
\begin{eqnarray}\label{P38}
r_{9}^{}\!\!&=\!\!&r_{9}'+r_{9}''~,\quad r_{9}'\,\succ\,r_{9}''~,
\end{eqnarray}
where
\begin{eqnarray}\label{P39}
\begin{array}{rcl}
r_{9}'\!\!&=\!\!&\displaystyle P_{+}\wedge\Bigl(\beta_{1}h+\alpha P_{2}
+\frac{\alpha\beta_{2}}{\beta_{1}}P_{1}\Bigr)\,+
\\[7pt]
&&+\,\displaystyle P_{1}\wedge\Bigl(\beta_{1}e_{+}+
\beta_{2}e'_{+}+\frac{\alpha\beta_{2}}{\beta_{1}}P_{+}\Bigr)~,
\\[12pt]
r_{9}'''\!\!&=\!\!&\displaystyle\chi P_{+}\wedge\Bigl(\beta_{1}e_{+}+
\beta_{2}e'_{+}+\frac{\alpha\beta_{2}}{\beta_{1}}P_{+}\Bigr)~.
\end{array}
\end{eqnarray}
The corresponding twisting two-tensor is given by the following formula
\begin{eqnarray}\label{P40}
F_{r_{9}^{}}\!\!&=\!\!&F_{r_{9}''}\,F_{r_{9}'}~,
\end{eqnarray}
where
\begin{eqnarray}\label{P41}
\begin{array}{rcl}
F_{r_{9}'}\!\!&=\!\!&\displaystyle\exp\Bigr(\bigr(\beta_{1}e_{+}^{}+
\beta_2e_{+}'+\frac{\alpha\beta_{2}}{\beta_{1}}P_{+}\bigr)\otimes P_{1}\Bigr)\,\times
\\[12pt]
&&\displaystyle\times\,\exp\Bigl(2\bigl(h+\frac{\alpha}{\beta_{1}}P_{2}
+\frac{\alpha\beta_{2}}{\beta_{1}^2}P_{1}\bigr)\otimes\sigma_{+}\Bigr)~,
\end{array}
\end{eqnarray}
\begin{eqnarray}\label{P42}
F_{r_{9}'''}\!\!&=\!\!&\displaystyle\exp\Bigl(\chi\bigl(e_{+}^{}+
\frac{\beta_{2}}{\beta_{1}}e_{+}'+\frac{\alpha\beta_{2}}{\beta_{1}^2}P_{+}\bigr)\wedge
\sigma_{+}\Bigr)~.\phantom{aaaaaaa}
\end{eqnarray}

10). The tenth $r$-matrix $r_{10}$,
\begin{eqnarray}\label{P43}
r_{10}^{}\!\!&=\!\!&\beta(P_{1}\wedge e'_{+}+P_{+}\wedge e_{+})+\alpha_{1}P_{-}\wedge
P_{1}+\alpha_{2} P_{+}\wedge P_{2}~,
\end{eqnarray}
can be presented as follows
\begin{eqnarray}\label{P44}
r_{10}^{}\!\!&=\!\!&r_{10}'+r_{10}''+r_{10}'''
\end{eqnarray}
where
\begin{eqnarray}\label{P45}
\begin{array}{rcl}
r_{10}'\!\!&=\!\!&\beta_{1}P_{+}\wedge e_{+}~,
\\[7pt]
r_{10}''\!\!&=\!\!&P_{2}\wedge(2\alpha_{1}P_{1}-\alpha_{2}P_{+})~,
\\[7pt]
r_{10}^{'''}\!\!&=\!\!&P_{1}\wedge(\beta_{1}e'_{+}-\alpha_{1}P_{-}+2\alpha_{1}P_{2})~.
\end{array}
\end{eqnarray}
All three $r$-matrices $r_{10}'$, $r_{10}''$, $r_{10}'''$ are Abelian and they have the
following subordination
\begin{eqnarray}\label{P46}
r_{10}'\,\succ\, r_{10}''~,\qquad r_{10}'+r_{10}''\,\succ\, r_{10}'''~.
\end{eqnarray}
The corresponding twisting two-tensor is given by the following formula
\begin{eqnarray}\label{P47}
F_{r_{10}^{}}\!\!&=\!\!&F_{r_{10}'''}\,F_{r_{10}''}\,F_{r_{10}'}~,
\end{eqnarray}
where
\begin{eqnarray}\label{P48}
\begin{array}{rcl}
F_{r_{10}^{'}}\!\!&=\!\!&\exp\bigr(\beta_{1}e_{+}\wedge P_{+}\bigr)~,
\\[9pt]
F_{r_{10}^{''}}\!\!&=\!\!&\exp\bigr((2\alpha_{1}P_{1}-\alpha_{2}P_{+})\wedge
P_{2}\bigr)~,
\\[9pt]
F_{r_{10}^{'''}}\!\!&=\!\!&\exp\bigl((\beta_{1}e'_{+}-\alpha_{1}P_{-}+2\alpha_{1}P_{2})
\wedge P_{1}\bigr)~.
\end{array}
\end{eqnarray}

11). The eleventh $r$-matrix $r_{11}$,
\begin{eqnarray}\label{P49}
r_{11}^{}\!\!&=\!\!&\beta_{1}P_{2}\wedge e_{+}+\alpha_{1}P_{+}\wedge P_{1}+\alpha_{2}
P_{-}\wedge P_{2}~,
\end{eqnarray}
is a sum of two subordinated $r$-matrices of Abelian type
\begin{eqnarray}\label{P50}
\begin{array}{rcl}
r_{11}^{}\!\!&=\!\!&r_{11}'+r_{11}''~,\quad r_{11}'\,\succ\, r_{11}''~,
\\[7pt]
r_{11}'\!\!&=\!\!&\alpha_{1}P_{+}\wedge P_{1}~,
\\[7pt]
r_{11}''\!\!&=\!\!&P_{2}\wedge(\beta_{1}e_{+}-\alpha_{2}P_{-})~.
\end{array}
\end{eqnarray}
The corresponding twisting two-tensor is given by the following formula
\begin{eqnarray}\label{P51}
F_{r_{11}^{}}\!\!&=\!\!&F_{r_{11}''}\,F_{r_{11}'}~,
\end{eqnarray}
where
\begin{eqnarray}\label{P52}
\begin{array}{rcl}
F_{r_{11}'}\!\!&=\!\!&\exp{(\alpha_{1}P_{1}\wedge P_{+})}~,
\\[9pt]
F_{r_{11}''}\!\!&=\!\!&\exp{\bigl((\beta_{1}e_{+}-\alpha_{2}P_{-})\wedge P_{2}\bigr)}~.
\end{array}
\end{eqnarray}

12). The two-tensor $\tilde{r}_{12}^{}$, which corresponds to the twelfth row of {\bf
Table 1}, has the form:
\begin{eqnarray}\label{P53}
\tilde{r}_{12}^{}\!\!&=\!\!&\beta_{1} P_{+}\wedge e_{+}+P_{-}\!\wedge(\alpha
P_{+}+\alpha_{1}P_{1}+\alpha_{2}P_{2})+ \tilde{\alpha} P_{+}\wedge P_{2}~.
\end{eqnarray}
We can show that $\tilde{r}_{12}^{}$ for $\alpha_{2}\neq0$ does not satisfy the equation
(\ref{P10}) that is the two-tensor $\tilde{r}_{12}^{}$ is not any classical $r$-matrix
provided that $\alpha_{2}\neq0$. In case of $\alpha_{2}=0$ the two-tensor is a classical
$r$-matrix and unfortunately in this case we had not success to construct the
corresponding twist.

13). The thirteen  $r$-matrix $r_{13}$,
\begin{eqnarray}\label{P54}
r_{13}^{}\!\!&=\!\!&\beta_{1}P_{0}\wedge h'+\alpha_{1}P_{0}\wedge P_{3}+\alpha_{2}
P_{1}\wedge P_{2}~,
\end{eqnarray}
is a sum of two subordinated $r$-matrices of Abelian type
\begin{eqnarray}\label{P55}
\begin{array}{rcl}
r_{13}^{}\!\!&=\!\!&r_{13}'+r_{13}''~,\quad r_{13}'\,\succ\, r_{13}''~,
\\[7pt]
r_{13}'\!\!&=\!\!&\alpha_{1}P_{0}\wedge P_{3}+\alpha_{2}P_{1}\wedge P_{2}~,
\\[7pt]
r_{13}''\!\!&=\!\!&\beta_{1}P_{0}\wedge h'~.
\end{array}
\end{eqnarray}
The corresponding twisting two-tensor is given by the following formula
\begin{eqnarray}\label{P56}
F_{r_{13}^{}}\!\!&=\!\!&F_{r_{13}''}\,F_{r_{13}'}~,
\end{eqnarray}
where
\begin{eqnarray}\label{P57}
\begin{array}{rcl}
F_{r_{13}'}\!\!&=\!\!&\exp{(\alpha_{1}P_{3}\wedge P_{0}+\alpha_{2}P_{2}\wedge P_{1})}~,
\\[10pt]
F_{r_{13}''}\!\!&=\!\!&\exp{(\beta_{1}h'\wedge P_{0})}~.
\end{array}
\end{eqnarray}

14). The fourteenth $r$-matrix $r_{14}$,
\begin{eqnarray}\label{P58}
r_{14}^{}\!\!&=\!\!&\beta_{1}P_{3}\wedge h'+\alpha_{1}P_{0}\wedge P_{3}+\alpha_{2}
P_{1}\wedge P_{2}~,
\end{eqnarray}
is a sum of two subordinated $r$-matrices of Abelian type
\begin{eqnarray}\label{P59}
\begin{array}{rcl}
r_{14}^{}\!\!&=\!\!&r_{14}'+r_{14}''~,\quad r_{14}'\,\succ\, r_{14}''~,
\\[7pt]
r_{14}'\!\!&=\!\!&\alpha_{1}P_{0}\wedge P_{3}+\alpha_{2}P_{1}\wedge P_{2}~,
\\[7pt]
r_{14}''\!\!&=\!\!&\beta_{1}P_{3}\wedge h'~.
\end{array}
\end{eqnarray}
The corresponding twisting two-tensor is given by the following formula
\begin{eqnarray}\label{P60}
F_{r_{14}^{}}\!\!&=\!\!&F_{r_{14}''}\,F_{r_{14}'}~,
\end{eqnarray}
where
\begin{eqnarray}\label{P61}
\begin{array}{rcl}
F_{r_{14}'}\!\!&=\!\!&\exp{(\alpha_{1}P_{3}\wedge P_{0}+\alpha_{2}P_{2}\wedge P_{1})}~,
\\[10pt]
F_{r_{14}''}\!\!&=\!\!&\exp{(\beta_{1}h'\wedge P_{3})}~.
\end{array}
\end{eqnarray}

15). The fifteenth  $r$-matrix $r_{15}$,
\begin{eqnarray}\label{P62}
r_{15}^{}\!\!&=\!\!&\beta_{1}P_{+}\wedge h'+\alpha_{1}P_{0}\wedge P_{3}+\alpha_{2}
P_{1}\wedge P_{2}~,
\end{eqnarray}
is a sum of two subordinated $r$-matrices of Abelian type
\begin{eqnarray}\label{P63}
\begin{array}{rcl}
r_{15}^{}\!\!&=\!\!&r_{15}'+r_{15}''~,\quad r_{15}'\,\succ\, r_{15}''~,
\\[7pt]
r_{15}'\!\!&=\!\!&\alpha_{1}P_{0}\wedge P_{3}+\alpha_{2} P_{1}\wedge P_{2}~,
\\[7pt]
r_{15}''\!\!&=\!\!&\beta_{1}P_{+}\wedge h'~.
\end{array}
\end{eqnarray}
The corresponding twisting two-tensor is given by the following formula
\begin{eqnarray}\label{P64}
F_{r_{15}^{}}\!\!&=\!\!&F_{r_{15}''}\,F_{r_{15}'}~,
\end{eqnarray}
where
\begin{eqnarray}\label{P65}
\begin{array}{rcl}
F_{r_{15}'}\!\!&=\!\!&\exp{(\alpha_{1}P_{3}\wedge P_{0}+\alpha_{2}P_{2}\wedge P_{1})}~,
\\[9pt]
F_{r_{15}''}\!\!&=\!\!&\exp{(\beta_{1} h'\wedge P_{+})}~.
\end{array}
\end{eqnarray}

16). The sixteenth  $r$-matrix $r_{16}$,
\begin{eqnarray}\label{P66}
r_{16}^{}\!\!&=\!\!&\beta_{1}P_{1}\wedge h+\alpha_{1}P_{0}\wedge
P_{3}+\alpha_{2}P_{1}\wedge P_{2}~,
\end{eqnarray}
is a sum of two subordinated $r$-matrices of Abelian type
\begin{eqnarray}\label{P67}
\begin{array}{rcl}
r_{16}^{}\!\!&=\!\!&r_{16}'+r_{16}''~,\quad r_{16}'\,\succ\, r_{16}''~,
\\[7pt]
r_{16}'\!\!&=\!\!&\alpha_{1}P_{0}\wedge P_{3}+\alpha_{2} P_{1}\wedge P_{2}~,
\\[7pt]
r_{16}''\!\!&=\!\!&\beta_{1}P_{1}\wedge h~.
\end{array}
\end{eqnarray}
The corresponding twisting two-tensor is given by the following formula
\begin{eqnarray}\label{P68}
F_{r_{16}^{}}\!\!&=\!\!&F_{r_{16}''}\,F_{r_{16}'}~,
\end{eqnarray}
where
\begin{eqnarray}\label{P69}
\begin{array}{rcl}
F_{r_{16}'}\!\!&=\!\!&\exp{(\alpha_{1}P_{3}\wedge P_{0}+\alpha_{2}P_{2}\wedge P_{1})}~,
\\[9pt]
F_{r_{15}''}\!\!&=\!\!&\exp{(\beta_{1}h\wedge P_{1})}~.
\end{array}
\end{eqnarray}

17). The seventeenth  $r$-matrix $r_{17}$,
\begin{eqnarray}\label{P70}
r_{17}^{}\!\!&=\!\!&\beta_{1} P_{+}\wedge h+\alpha P_{1}\wedge
P_{2}+\alpha_{1}P_{+}\wedge P_{1}~,
\end{eqnarray}
is a sum of two subordinated $r$-matrices of Jordanian and Abelian types
\begin{eqnarray}\label{P71}
\begin{array}{rcl}
r_{17}^{}\!\!&=\!\!&r_{17}'+r_{17}''~,\quad r_{17}'\,\succ\, r_{17}''\,\prec\,r_{17}'~,
\\[7pt]
r_{17}'\!\!&=\!\!&P_{+}\wedge(\beta_{1}h+\alpha_{2}P_{1})~,
\\[7pt]
r_{17}''\!\!&=\!\!&\alpha_{1} P_{1}\wedge P_{2}~.
\end{array}
\end{eqnarray}
The corresponding twisting two-tensor is given by the following formula
\begin{eqnarray}\label{P72}
F_{r_{17}^{}}\!\!&=\!\!&F_{r_{17}''}\,F_{r_{17}'}\,=\,F_{r_{17}'}\,F_{r_{17}''}~,
\end{eqnarray}
where
\begin{eqnarray}\label{P73}
\begin{array}{rcl}
F_{r_{17}'}\!\!&=\!\!&\displaystyle\exp{\Bigl(\bigl(h+
\frac{\alpha_{2}}{\beta_{1}}P_{1}\bigr)\otimes\sigma_{+}\Bigr)}~,
\\[12pt]
F_{r_{17}''}\!\!&=\!\!&\exp{(\alpha_{1} P_{2}\wedge P_{1})}~.
\end{array}
\end{eqnarray}

18). The eighteenth $r$-matrix $r_{18}$,
\begin{eqnarray}\label{P74}
r_{18}^{}\!\!&=\!\!&P_{+}\wedge(\beta_{1}h+\beta_{2}h')+\alpha P_{1}\wedge P_{2}~,
\end{eqnarray}
is a sum of two subordinated $r$-matrices of Abelian and Jordanian types
\begin{eqnarray}\label{P75}
\begin{array}{rcl}
r_{18}^{}\!\!&=\!\!&r_{18}'+r_{18}''~,\quad r_{18}'\,\succ\, r_{18}''~,
\\[7pt]
r_{18}'\!\!&=\!\!&\alpha P_{1}\wedge P_{2}~,
\\[7pt]
r_{18}''\!\!&=\!\!&P_{+}\wedge(\beta_{1}h+\beta_{2} h')~,
\end{array}
\end{eqnarray}
The corresponding twisting two-tensor is given by the following formula
\begin{eqnarray}\label{P76}
F_{r_{18}^{}}\!\!&=\!\!&F_{r_{18}''}\,F_{r_{18}'}~,
\end{eqnarray}
where
\begin{eqnarray}\label{P77}
\begin{array}{rcl}
F_{r_{18}'}\!\!&=\!\!&\exp{(\alpha P_{2}\wedge P_{1})}~,
\\[10pt]
F_{r_{18}''}\!\!&=\!\!&\displaystyle\exp{\Bigl(\bigl(h+
\frac{\beta_{2}}{\beta_{1}}h'\bigr)\otimes\sigma_{+}\Bigr)}~.
\end{array}
\end{eqnarray}

19). -- 21). The $r$-matrices:
\begin{eqnarray}\label{P78}
r_{19}^{}\!\!&=\!\!&\alpha P_{1}\wedge P_{+}~,
\\[7pt]\label{P79}
r_{20}^{}\!\!&=\!\!&\alpha P_{1}\wedge P_{2}~,
\\[7pt]\label{P80}
r_{21}^{}\!\!&=\!\!&\alpha_{1} P_{0}\wedge P_{3}+\alpha_{2}P_{1}\wedge P_{2}
\end{eqnarray}
are Abelian and their corresponding twists are given by the simple exponential formulas
\begin{eqnarray}\label{P81}
F_{r_{19}^{}}\!\!&=\!\!&\exp{(\alpha P_{+}\wedge P_{1})}~,
\\[9pt]\label{P82}
F_{r_{20}^{}}\!\!&=\!\!&\exp{(\alpha P_{2}\wedge P_{1})}~,
\\[9pt]\label{P83}
F_{r_{21}^{}}\!\!&=\!\!&\exp{(\alpha_{1} P_{3}\wedge P_{0}+\alpha_{2}P_{2}\wedge
P_{1})}~.
\end{eqnarray}

\section*{Conclusion}
In this paper twists which describe multiparameter quantum deformations of the
Poincar\'{e} algebra were obtained in explicit form. These twists correspond to eighteen
out of the twenty classical $r$-matrices of the Zakrzewski's classification, which
satisfy the homogeneous classical Yang-Baxter equation. For this aid we used the notation
of subordination for the classical $r$-matrices. These results can be extend on the
Poincar\'{e} superalgebra. It should be noted that up to now there is not any
classification of the classical $r$-matrices in spirit of the classification by S.
Zakrzewski. However it turns out that it is possible to extend the Zakrzewski's
classification to the Poincare superalgebra by an addition of supercharges terms to the
classical $r$-matrices $r_i$ ($i = 1,2,\ldots,21$). Theses results and corresponding
twisting functions will be presented in a future paper.


\end{document}